\documentclass[final,3p,times]{elsarticle}
\usepackage{amssymb}
\usepackage{algorithm}
\usepackage{algpseudocode}
\usepackage{amsmath}
\usepackage{xurl}
\usepackage{booktabs}

\usepackage[english]{babel}
\usepackage{amsthm}

\usepackage{gensymb}

\journal{Computer Methods in Applied Mechanics and Engineering}

\begin{document}

\begin{frontmatter}

\title{Biomimetic IGA neuron growth modeling with neurite morphometric features and CNN-based prediction}

\author[a]{Kuanren Qian}
\author[a]{Ashlee S. Liao}
\author[b]{Shixuan Gu}
\author[a,b,c]{Victoria A. Webster-Wood}
\author[a,b]{Yongjie Jessica Zhang}

\address[a]{Department of Mechanical Engineering, Carnegie Mellon University, 5000 Forbes Ave, Pittsburgh, PA 15213, USA}
\address[b]{Department of Biomedical Engineering, Carnegie Mellon University, 5000 Forbes Ave, Pittsburgh, PA 15213, USA}
\address[c]{McGowan Institute for Regenerative Medicine, University of Pittsburgh, 450 Technology Drive, Pittsburgh, PA 15219, USA}

\begin{abstract}
Neuron growth is a complex, multi-stage process that neurons undergo to develop sophisticated morphologies and interwoven neurite networks. Recent experimental research advances have enabled us to examine the effects of various neuron growth factors and seek potential causes for neurodegenerative diseases, such as Alzheimer's disease, Parkinson's disease, and amyotrophic lateral sclerosis. A computational tool that studies the neuron growth process could shed crucial insights on the effects of various factors and potentially help find a cure for neurodegeneration. However, there is a lack of computational tools to accurately and realistically simulate the neuron growth process within reasonable time frames. Bio-phenomenon-based models ignore potential neuron growth factors and cannot generate realistic results, and bio-physics-based models require extensive, high-order governing equations that are computationally expensive. In this paper, we incorporate experimental neurite features into a phase field method-based neuron growth model using an isogeometric analysis collocation (IGA-C) approach. Based on a semi-automated quantitative analysis of neurite morphology, we obtain relative turning angle, average tortuosity, neurite endpoints, average segment length, and the total length of neurites. We use the total neurite length to determine the evolving days \textit{in vitro} (DIV) and select corresponding neurite features to drive and constrain the neuron growth. This approach archives biomimetic neuron growth patterns with automatic growth stage transitions by incorporating corresponding DIV neurite morphometric data based on the total neurite length of the evolving neurite morphology.
Furthermore, we built a convolutional neural network (CNN) to significantly reduce associated computational costs for predicting complex neurite growth patterns. Our CNN model adopts a customized convolutional autoencoder as the backbone that takes neuron growth simulation initializations and target iteration as the input and predicts the corresponding neurite patterns. This approach achieves high prediction accuracy  (97.77\%) while taking 7 orders of magnitude less computational times when compared with our IGA-C neuron growth solver.
\end{abstract}


\end{frontmatter}

\section{Introduction}
In recent years, we have witnessed significant outgrowth in the neuron growth research field, and multiple attempts to prevent neurodegenerative with neuron growth factors have been made. Understanding nervous system development is critical in searching for potential neurodegenerative disease treatments. Recent molecular biology research has revealed many possible effects of neuron growth factors during disease pathogenesis \cite{dugger_pathology_2017,brown_neurodegenerative_2005}. Alzheimer's disease, Parkinson's disease, and amyotrophic lateral sclerosis are three devastating neurodegenerative disorders with high morbidity and mortality for patients \cite{checkoway_neurodegenerative_2011}. These diseases arise when neurons progressively lose connections or functions due to alterations of certain neuron growth factors \cite{berg_new_1984,connor_role_1998}. Because most mammalian neurons cannot intrinsically regenerate, they can not repair or replace themselves, making traditional treatments ineffective against neurodegeneration \cite{reimer_motor_2008,steward_neural_2012}. Although evidence shows that certain factors from developing neurons can help protect mature neurons from degeneration \cite{connor_role_1998}, there is a lack of study of the specific functional role of these factors \cite{elliott_motor_1996}. Therefore, understanding the factors involved during neuron growth is vital for seeking potential neurodegeneration treatments. Neuron growth is a complex multi-stage process consisting of multiple stages of development, each with unique growth behaviors \cite{van_ooyen_modeling_2003}. Primary rat hippocampal neurons \textit{in vitro} have exhibited five key developmental stages in their morphology \cite{dotti_establishment_1988}. In Stage 1, lamellipodia protrude from the initial neuron cell and extend into neurites of similar length in Stage 2. The longest neurite will differentiate into a longer axon in Stage 3, while other neurites slowly grow in Stage 4. Stage 5 is neuron maturation, in which neurons grow intricate patterns and form complex neurite networks. This process spans several days \textit{in vitro} (DIV), and each stage exhibits drastically different growth behaviors driven by intracellular and extracellular biophysics processes that involve a wide range of neuron growth factors \cite{hentschel_instabilities_1998,liao_quantitative_2022}.

Considering the aggregating costs and arduous procedures needed for extensive cell culture experiments, mathematical modeling of early neuron growth stages has been proposed to study neuron growth factor during the initial neurite outgrowth \cite{hentschel_instabilities_1994}, axon differentiation \cite{krottje_mathematical_2007,pearson_mathematical_2011}, and axon guidance \cite{aeschlimann_biophysical_2001}. Among different neuron growth models, there exist two major schools of thought. One approach models neuron growth using stochastic methods that follow phenomenological results \cite{eberhard_neugen_2006,van_ooyen_independently_2014}. Along this direction, some literature models the growth process using axon steering theory based on filopodia \cite{goodhill_predicting_2004}, external repulsive cues \cite{maskery_growth_2004}, and stochastic mechanism \cite{koene_netmorph_2009}. Generalized neurite characteristics are utilized based on morphology \cite{cuntz_one_2010,donohue_comparative_2008}, and the surrounding substrate is incorporated \cite{torben-nielsen_context-aware_2014}. While efficient to compute, these methods take limited biophysics into account. In contrast, models that attempt to incorporate the effects of different neuron growth factors based on biophysical mechanisms behind neurite elongations are prone to high computational cost and numerical instability \cite{otoole_physical_2008,graham_mathematical_2006}. Despite the aforementioned approaches to model the neuron growth process, most of them lack the ability to efficiently and accurately model realistic neuron growth. Many rely heavily on broad assumptions or are associated with expensive computational costs, yet are still unable to capture complex neurite morphology changes throughout the growth process.

Inspired by \cite{takaki_phase_field_2015}, we developed a phase field method-based neuron growth model with distinctive growth stage transitions \cite{qian_modeling_2022} to address these limitations in neuron growth modeling. Our isogeometric analysis collocation (IGA-C) neuron growth model archives multi-stage neuron growth by considering the effect of intracellular tubulin concentration. By adjusting the assembly and disassembly rate of the tubulin transport coupled with the phase field model, the model could simulate axon differentiation through different neurite elongation rates based on tubulin concentration. Yet, despite our implementation of biophysics-based approaches, our neuron growth model relies on arbitrarily set growth stage transitions and can not capture similar biomimetic growth behaviors without extensive manual parameter adjustments for a particular type of neuron. As a follow-up, in this paper we incorporate experimentally observed neurite morphometric features into the biophysics-based IGA-C neuron growth model to simulate biomimetic neuron growth behaviors with intrinsic growth stage transitions determined by neurite morphology. Based on the semi-automated quantitative neurite morphometric evaluation \cite{liao_semi-automated_2022}, the developed neuron growth model uses the total neurite length to determine the evolving DIV and select neurite morphometrics of the corresponding DIV to constrain and drive neurite growth behaviors. The model can be easily adapted to different types of neurons, given the corresponding neurite morphometrics.

Note that our neuron growth simulations with incorporated neurite features require days to complete, and the associated computational costs arise rapidly for multiple neuron cases as the domain increases in size. In addition, the neurite morphometric feature implementation requires evaluating each incorporated neurite feature within each iteration. In recent years, deep learning has been used in solving nonlinear high-dimensional partial differential equations (PDEs) \cite{han_solving_2018}, and physics-informed neural networks \cite{raissi_physics-informed_2019} were developed. They have significantly alleviated computational costs without undermining analysis accuracy when studying a 2D reaction-diffusion system \cite{li_reaction_2020}, intracellular material transport process \cite{li_deep_2021}, and traffic jams in complex 3D neuron structures \cite{li2023isogeometric}. In this paper, we propose to use a convolutional neural network-based (CNN) surrogate model to accurately predict intricate neurite growth patterns with significantly lower computational costs than conventional simulations. The main contributions of the presented neurite morphometric feature-driven IGA-C growth model and CNN-based prediction include:
\begin{itemize}    
    \item Incorporating experimentally observed neurite features into the IGA-C phase field neuron growth model to drive and constrain biomimetic neuron growth behaviors with seamless intrinsic growth stage DIV transitions;
    \item An IGA-C neuron growth model that simulates rat hippocampal neuron growth process and is adaptable to different types of neurons given neurite morphometric data of the growth behaviors; and
    \item A deep learning model based on CNN architecture, capable of accurate and fast predictions of complex neurite networks. The model predicts neurite growth across different growth stages without an iterative process.
\end{itemize}
The remainder of this paper is organized as follows. Section 2 reviews IGA-C-based phase field neuron growth. Section 3 elaborates on the neurite morphometric features-driven neuron growth model that simulates biomimetic neurite growth with intrinsic growth stage transition. Section 4 presents a CNN-based surrogate model for neuron growth and showcases prediction results. Finally, we draw conclusions and discuss future work.

\section{Review of IGA-C-based phase field neuron growth modeling}

In this section, we review the phase field neuron growth model, IGA, and collocation method. IGA and phase field methods are powerful numerical techniques utilized in simulating complex engineering problems. IGA is a high-order numerical technique used in computational mechanics to approximate governing equations continuously, eliminating the need for discretization in traditional finite element methods \cite{hughes_isogeometric_2005}. The phase field method solves evolving boundary problems like fracture and dendrite solidification \cite{takaki2014phase}. Combining these two techniques can effectively simulate neurite growth with high accuracy. 

\subsection{Phase field neuron growth model}

In our previous work, we developed an IGA-C-based phase field model to simulate the multi-stage neuron growth process by introducing the effect of tubulin concentration \cite{qian_modeling_2022}. By considering intracellular transport during the growth, our model differentiates the longest neurite into an axon. We simulated different growth behaviors during each stage to reproduce growth cone behaviors at neurite tips. By treating the neuron growth domain as a binary phase field $\phi$ with phase ``1'' representing the neuron and phase ``0'' representing the surrounding medium, we solved the phase field governing equations using the IGA-C method to simulate neurite growth patterns based on given neuron configurations and parameters. In particular, we used IGA-C to solve the phase field equation (Eq.~\ref{phase field equation}) coupled with the tubulin equation (Eq.~\ref{tubulin equation}) and the temperature equation (Eq.~\ref{temperature equation}) through an energy term $E$ (Eq.~\ref{energy equation}) and competitive tubulin consumption term $\frac{dL}{dt}$ (Eq.~\ref{dLdt equation}). We provide a detailed list of variables in Table~\ref{Table: phase field parameter} for clarity. The main governing equations are defined as:
\begin{align}
&\frac{\partial \phi}{\partial t} = M_{\phi} [\bigtriangledown\cdot(a(\Psi)^2\bigtriangledown \phi) - \frac{\partial} {\partial{x}}(a(\Psi) \frac{\partial a(\Psi)}{\partial \Psi} \frac{\partial \phi}{\partial y}) + \frac{\partial} {\partial{y}}(a(\Psi) \frac{\partial a(\Psi)}{\partial \Psi} \frac{\partial \phi}{\partial x}) + \phi(1-\phi)(\phi - 0.5 + E + 6H |\bigtriangledown \theta|)], \label{phase field equation}\\
&\frac{\partial (\phi \,c_{tub})}{\partial t} =  \delta_t \bigtriangledown\cdot \: (\phi \, \bigtriangledown c_{tub}) - \mathbf{\alpha}_{t} \cdot \bigtriangledown (\phi \, c_{tub}) - \beta_{t} (\phi \, c_{tub}) + \epsilon_0 \frac{|\bigtriangledown(\phi_0)|^2}{\int{|\bigtriangledown(\phi_0)|^2} \, d \Omega}, \label{tubulin equation}\\
&\frac{\partial T}{\partial t} = \bigtriangledown^2T + K\frac{\partial \phi}{\partial t}, \label{temperature equation}\\
&E =  \frac{\alpha}{\pi}	\tan^{-1}(H_\epsilon(\frac{d L}{d t}) \gamma \bigtriangleup T),\label{energy equation}\\
&\frac{dL}{dt} = r_{g} \, c_{tub} - s_{g}, \label{dLdt equation}
\end{align} where $M_{\phi}$ is the phase field mobility coefficient, $a(\Psi)$ is the anisotropy coefficient for the gradient, $E$ is the energy term for phase field growth, $H$ is a constant value, and $\theta$ is the orientation term. In Eq. \ref{tubulin equation}, $\delta_t$ is the tubulin diffusion rate, $\alpha_t$ is the tubulin active transport coefficient, $\beta_t$ is the tubulin decay coefficient, and $\epsilon_0 \frac{|\bigtriangledown(\phi_0)|^2}{\int{|\bigtriangledown(\phi_0)|^2} \, d \Omega}$ is the constant tubulin production source term, where $\phi_0$ is the initial phase field variable and $\epsilon_0$ is the dimensionless production coefficient term. In Eq. \ref{temperature equation}, $\frac{\partial T}{\partial t}$ is the undercooling temperature equation, and $K$ is the dimensionless latent heat coefficient. In Eq.~\ref{energy equation}, $E$ is the energy term responsible for driving cell growth at the interface between phases. $\frac{\alpha}{\pi}$ is a scaling coefficient, $H_\epsilon$ is the Heaviside step function, and $\gamma$ is the phase field interfacial energy constant. In Eq.~\ref{dLdt equation}, $r_g$ and $s_g$ are the assembly and disassembly rate of tubulin concentration within the cell boundary \cite{van_ooyen_competition_2001} to incorporate the effect of intracellular tubulin concentration.

\begin{table}[ht]
\caption{Parameters utilized in the phase field neuron growth model.}
\vspace{-0.7cm}
\label{Table: phase field parameter}
\begin{center}
\resizebox{\columnwidth}{!}{\begin{tabular}{c l c|c l c}
\toprule
\textbf{Parameter} & \textbf{Description} & \textbf{Value} &
\textbf{Parameter} & \textbf{Description} & \textbf{Value}\\
\hline
$\phi$ & phase field variable & - & $\gamma$ & phase field interface energy constant & $10$ \\  
$M_{\phi}$ & Mobility coefficient & $60$ & $c_{tub}$ & Tubulin concentration & - \\
$a(\Psi)$ & Neurite morphology gradient coefficient & - & $\delta_t$ & Tubulin diffusion rate & $4$ ($\mu m^2/h$) \\
$H$ & Orientation constant coefficient & $0.007$ & $\alpha_t$ & Tubulin active transport rate & $0.001$ ($\mu m/h$) \\
$\theta$ & Neurite growth orientation angle & $[0,1]$ & $\beta_t$ & Tubulin decay coefficient & $0.001$ ($1/h$) \\
$T$ & Temperature & - & $\epsilon_0$ & Tubulin production coefficient & $15$ \\
$K$ & Dimensionless latent heat & $2$ & $\phi_0$ & Initial phase field variable & - \\
$E$ & Energy & - & $\frac{dL}{dt}$ & Competitive tubulin consumption & - \\  
$\frac{\alpha}{\pi}$ & Scaling coefficient & $0.2865$ & $r_g$ & Tubulin assembly rate & $5$ \\  
$H_{\epsilon}$ & Heaviside function & - & $s_g$ & Tubulin disassembly rate & $0.1$ \\  
\bottomrule
\\[-1em] 
\multicolumn{6}{l}{Note: Default value is given if the parameter requires initialization. For dimensionless parameters, we show the default value without units.}
\end{tabular}}
\end{center}
\vspace{-0.75cm}
\end{table}

The neuron growth model achieves neurite outgrowth by iteratively solving interface evolution using the phase field governing equation (Eq.~\ref{phase field equation}). The first few terms involving $a(\Psi)$ introduce anisotropy into the neuron growth model \cite{takaki_phase_field_2015}. The following term involving $E$ is a double well function with a coupling orientation term, $6H |\bigtriangledown \theta|$, introduced in \cite{ren_controllable_2018} to break dendrite growth symmetry. In the phase field model, $\theta$ indicates the change in the direction of the extending neurites. The initial orientation $\theta$ is set as random values between $[0, 1]$ throughout the domain and the cell is initialized as a filled circle at the center of the binary domain $\phi$, where we consider $\phi_0 = 1$ in the cell and $\phi_0 = 0$ in the medium. The seed radius of the cell is set as $r_0 = 20\,\triangle x$. We initialize the tubulin concentration in the cell as $c_{tub} = \frac{1}{2}(1+ tanh((r_0-r)/2))$. In our earlier neuron growth model \cite{qian_modeling_2022}, we
\begin{itemize}    
    \item Developed a phase field model coupled with intracellular tubulin concentration to simulate the first four stages of neuron growth, including lamellipodia formation, initial neurite outgrowth, axon differentiation, and dendrite formation; and
    \item Implemented relative turning angle from experimental observations to guide neurite growth and demonstrated similar end-stage neuron morphology reproduction; and
    \item Extended the model to simulate the growth of neural circuits, where multiple neurons build neurite interactions and form neural networks.
\end{itemize} 

\vspace{0.2cm}
\textbf{Discussion}. Despite the attempts made in the aforementioned model \cite{qian_modeling_2022} to simulate growth behaviors, the model still relies on arbitrarily-set iterations for each growth stage and can not accurately depict intrinsic growth stage transition. In addition, the model develops excessive branchings at later stages not observed in experimental cultures. Moreover, parameters used in the model were manually tuned for rat hippocampal neuron behaviors, making it less applicable to other types of neurons. These limitations undermine the practical implementations of the model. In this paper, we propose incorporating neurite morphometric features from experimental cultures as a modular component into the neuron growth model to overcome these limitations. Our proposed approach will improve the model automation and its ability to simulate more biomimetic growth behaviors of a diverse range of neurons.

\subsection{Isogeometric-collocation (IGA-C) method}
In this work, we utilize the IGA-C method to solve the phase field neuron growth model for its flexibility and accuracy \cite{gomez_accurate_2014,schillinger_isogeometric_2015}. High-fidelity modeling for complex neuron morphologies is a significant challenge in computational biology \cite{zhang_challenges_2013,zhang_geometric_2018}. Using high-order spline basis functions, we can generate an accurate analysis of the phase field model. Together with the collocation method, we directly solve the strong form of the phase field equation on Greville Abscissae collocation points \cite{auricchio_isogeometric_2010,casquero_isogeometric_2016}.

\textbf{B-splines basis function.} In this paper, we choose cubic B-splines as the basis functions for our IGA-C solver \cite{piegl1996nurbs,wei2017truncated}. For a univariate B-spline curve of degree $p$, it is defined using a sequence of nondecreasing real numbers, $u_i$, to construct an open knot vector $U = \{u_1, u_2, \cdots u_{n+p+1}\}$, in which $p$ is the B-spline order, and $n$ is the number of basis functions. The basis function $N_{i,p}(u)$ is defined based on the knot vectors, we have: 
\begin{equation}
    N_{i,0}(u) = 
    \begin{cases}
        1, & \text{if $u_i \le u \le u_{i+1}$} \\
        0, & \text{otherwise,}
    \end{cases}
\end{equation}
\begin{equation}
    N_{i,p}(u) = \frac{u-u_i}{u_{i+p}-u_i}N_{i,p-1}(u) + \frac{u_{i+p+1}-u}{u_{i+p+1}-u_{i+1}}N_{i+1,p-1}(u), \quad p=1,2,...
\end{equation} where $N_{i,0}(u)$ is piece-wise linear constant. $N_{i,p}(u)$ is recursively defined based on a combination of preceding basis functions. With control points $P = \{P_{i}\}_{i=1}^n$, we can construct a $p^{th}$-order B-spline curve $C(u)$ as:
\begin{equation}
    C(u) = \sum^n_{i=0}{N_{i,p}(u)P_i}, \quad 0 \le u \le 1 
\end{equation} where $N_{i,p}(u)$ is the B-spline basis function of $p$ degree, defined over the knot vector $U$.
Given an arbitrary $u$ value in the knot vector $U$, we can compute the point on the B-spline curve by multiplying the value of every non-zero basis function $N_{i,p}(u)$ with its control point $P_i$ and taking a summation. For a 2D surface, we define the basis functions as the tensor product of two univariate B-splines. 

\textbf{Greville abscissae collocation points.} The Greville abscissae collocation method is a numerical technique for solving strong form of PDEs \cite{jia_adaptive_2019,casquero_isogeometric_2016}.  Using Greville abscissae points, it is relatively straightforward to implement and accurate for problems with smooth solutions. It is also adaptable to a wide range of high-order differential equations. For IGA, collocation methods have been recently shown to perform well as an alternative to Galerkin methods \cite{anitescu_isogeometric_2015, casquero_hybrid_2016} while maintaining higher order convergence rates \cite{johnson2005higher}. The Greville abscissae point $\hat{u_i}$ is defined based on knots in the knot vector $U = \{u_1, u_2, \cdots u_{n+p+1}\}$ as:
\begin{equation}
    \hat{u_i} = \frac{u_{i+1}+...+u_{i+p}}{p}, \quad i \in{[1,n]}
\end{equation} which can be directly computed on given B-spline surfaces \cite{johnson2005higher}. In the context of neuron growth, we can rewrite each collocation point as $\hat{\phi} = \{\hat{\phi}_u,\hat{\phi}_v\}$:
\begin{equation}
    \hat{\phi}_u = \frac{\sum_{i+1}^{i+p}u}{p} \; \text{ and } \; \hat{\phi}_v = \frac{\sum_{j+1}^{j+p}v}{p},
\end{equation} where $\hat{\phi}_u$ and $\hat{\phi}_v$ are the components along each parametric direction of the collocation point $\hat{\phi}$. By directly solving the strong form of the phase field neuron growth model, we obtain accurate and smooth results. 

\section{Neurite morphometric features-driven neuron growth}

In this section, we first review experimental neuron culture procedures and semi-automated neurite morphometric evaluation. Then, we present the computational pipeline for the feature-driven neuron growth model and go through each implemented module in detail. Finally, we discuss our simulation results.

\subsection{Experimental neuron culture and semi-automated neurite morphometric evaluation}

\begin{figure}[ht]
\centering
\includegraphics[width=\linewidth]{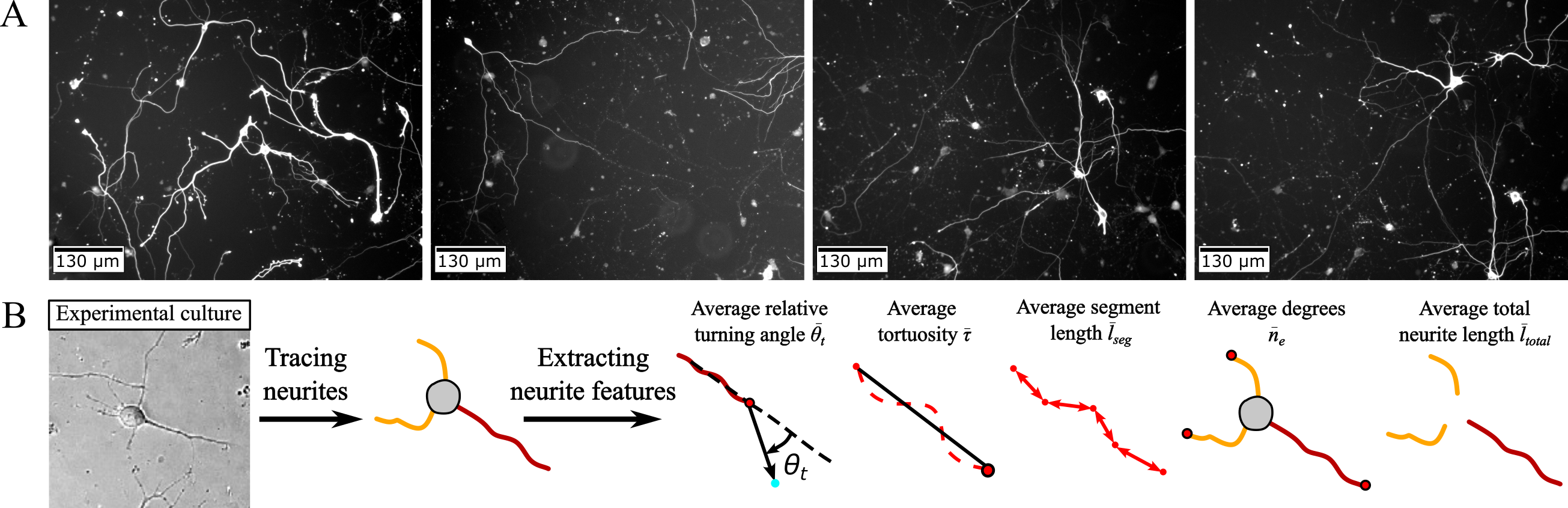}
\caption{Neuron growth images and neurite morphometric features. (A) Fluorescent images of rat hippocampal neurons to illustrate the complexity of neuron structure. (B) Neurite morphometric features were extracted from experimental images of neuron cultures using semi-automated neurite morphometric evaluation \cite{liao_semi-automated_2022}, including average relative turning angle $\bar \theta_t$, average tortuosity $\bar \tau$, average segment length $\bar l_{seg}$, average degree $\bar n_e$, and average total neurite length $\bar l_{total}$.}
\label{figure: Neurite_morphometric_features}
\vspace{-0.3cm}
\end{figure}

\textbf{Neuron cell culture}. Neuron growth exhibits complex morphological characteristics as illustrated in Figure~\ref{figure: Neurite_morphometric_features}A, and behavior changes throughout the process. For the work presented here, data was used from an existing experimental neuron growth dataset of bright-field microscopy images \cite{liao_semi-automated_2022}. Briefly, embryonic rat hippocampal neurons were cultured over 6 DIV, and the first week of \textit{in vitro} morphological development was imaged and analyzed. Cryopreserved primary, embryonic-day 18 (E18) rat hippocampal neurons (A36513, Gibco, USA) were cultured in dishes coated with poly-D-lysine (P6407, Sigma-Aldrich, USA) following manufacturer protocol \cite{2018B-27System}. Neurons at a density of $10,000$ cells/$cm^2$ were seeded in Neurobasal Plus (A3582901, Gibco, USA) supplemented with 2\% B-27 Plus (A3582801, Gibco, USA). Except during the media change and microscopy imaging, the incubation was conducted at 37 $^{\circ}C$ with $5\%$ $CO_2$. $50\%$ culture media was replenished with fresh media after 24 hours. The imaging process utilized the Echo Revolve Microscope (Echo Revolve $|$ R4, inverted, BICO - The Bio Convergence Company, USA) with a 12-megapixel color camera in the inverted configuration. The magnification was set at 20X or 40X, and cells were imaged at 0.5, 1, 2, 3, 4, and 6 DIV. 

\textbf{Experimentally-derived neurite morphometric features}. Neurite morphometric features used to drive the growth model in this work were taken from our prior neuron growth dataset \cite{liao_semi-automated_2022}. We leveraged the distributions of representative neurite morphological features per cell (Figure~\ref{figure: Neurite_morphometric_features}B), including average relative turning angle $\bar \theta_t$,  average tortuosity $\bar \tau$, average segment length $\bar l_{seg}$, neurite degrees (the number of neurite endpoints/tips $\bar n_e$), and the average total length $\bar l_{total}$ as a dataset shown in Table \ref{table: neurite_feature_data} \cite{liao_quantitative_2022,liao_semi-automated_2022}. The relative turning angles $\bar \theta_t$ is the angle change between consecutive change point segments. The tortuosity ${\bar \tau}$ is calculated by dividing the length of the neurite tracing by the distance between its endpoints. The average segment length $\bar l_{seg}$ is the average distance between change points of all neurites for a given cell. Degree $\bar n_e$ is the number of neurite endpoints for a given cell. The total length $\bar l_{total}$ is the length of all traced neurites part of a given cell.

\begin{table}[ht]
\centering
\caption{\label{table: neurite_feature_data} Neurite morphometric features statistics from DIV 0.5 to 6 from rat hippocampal neurons \cite{liao_semi-automated_2022}.}
\begin{tabular}{c|c c|c c|c c|c c|c}
\toprule
DIV &  \multicolumn{2}{|c|}{$\bar \theta_t (\degree) $} & \multicolumn{2}{|c|}{$\bar \tau$} & \multicolumn{2}{|c|}{$\bar n_e$} & \multicolumn{2}{|c|}{$\bar l_{seg} (\mu m) $} & $\bar l_{total} (\mu m)$ \\
\cline{2-9}
{}  & $\mu$ & $\sigma$ & $Q_1$ & $Q_3$ & $Q_1$ & $Q_3$ & $Q_1$ & $Q_3$ 
 & {}\\
\hline
$0.5$ & $22.18$ & $11.65$ & $1.0225$ & $1.0776$ & $1$ & $2.5$ & $5.05$ & $7.82$ & $27.53$\\
\hline
$1$ & $22.82$ & $10.82$ & $1.0161$ & $1.0757$ & $1$ & $3$ & $5.58$ & $8.36$ & $36.54$\\
\hline
$1.5$ & $21.34$ & $9.69$ & $1.0254$ & $1.0507$ & $1$ & $4$ & $5.64$ & $9.39$ & $53.19$\\
\hline
$2$ & $22.32$ & $8.01$ & $1.0283$ & $1.0685$ & $1$ & $4$ & $6.39$ & $10.35$ & $84.34$\\
\hline
$3$ & $22.88$ & $6.78$ & $1.0300$ & $1.0725$ & $2$ & $6$ & $8.35$ & $10.95$ & $155.13$\\
\hline
$4$ & $21.28$ & $4.91$ & $1.0341$ & $1.0623$ & $3$ & $7$ & $6.65$ & $10.86$ & $218.74$\\
\hline
$6$ & $20.32$ & $2.65$ & $1.0302$ & $1.0498$ & $6$ & $10$ & $10.24$ & $12.36$ & $554.73$\\
\bottomrule
\multicolumn{10}{l}{Note: $\mu$ and $\sigma$ are the mean and standard deviation, respectively. $Q_1$ and $Q_3$ are the } \\
\multicolumn{10}{l}{1st and 3rd quartiles of the features distribution.} \\
\end{tabular}
\end{table}

\subsection{Incorporating neurite morphometric features}
\begin{figure}[htbp]
\centering
\includegraphics[width=\linewidth]{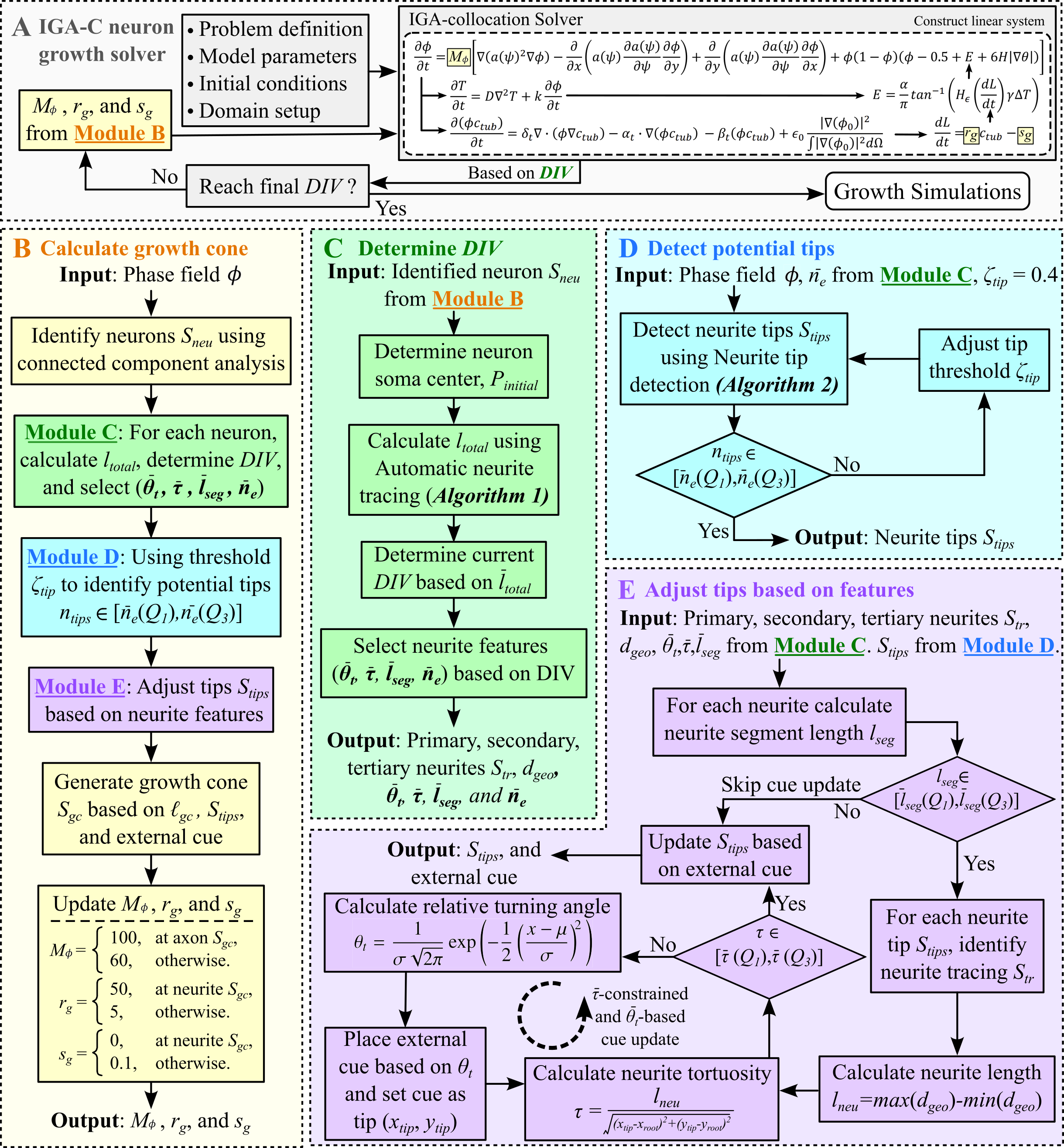}
\caption{IGA-C neuron growth computational pipeline incorporating neurite morphometric features extracted from experiments. (A) IGA-C phase field neuron growth solver. (B) Updating neuron growth cone to guide neurite outgrowth. (C) Intrinsic growth stage transition based on Algorithm~\ref{Algorithm: Neurite_tracing_algorithm} and determined DIV. (D) Detecting potential tips using Algorithm~\ref{Algorithm: tip_detection}. (E) Neurite morphometric features-based tip adjustments.}
\label{figure: Computational_pipeline}
\vspace{-0.5cm}
\end{figure}

Our previous neuron growth model \cite{qian_modeling_2022} only utilizes relative turning angles $\theta_t$ extracted from experimental data without adjustments throughout the entire simulation process. Therefore, it can not intrinsically transition between adjacent neuron growth stages (DIV), leading to unnecessarily dense neurite branching during the simulation. In this paper, we incorporate DIV-based neurite morphometric features from embryonic rat hippocampal neuron cultures into our neuron growth model to better capture neurite morphological characteristics that vary across multiple growth stages intrinsically. We provide a detailed explanation of parameters used in the feature-driven neuron growth procedure in Table~\ref{feature-driven parameter table}, and show the overall computational pipeline in Figure~\ref{figure: Computational_pipeline}. Our feature-driven neuron growth consists of five modules: IGA-C neuron growth solver (Module A), growth cone calculation (Module B), DIV determination (Module C), potential tip detection (Module D), and feature-based tip adjustment (Module E). Module A has been reviewed in Section 2.1. We explain Modules B-E in detail below. Parameters in Table 3 are categorized based on these four modules.

\begin{table}[ht]
\caption{Parameters utilized in feature-driven neuron growth procedure.}
\vspace{-0.7cm}
\label{feature-driven parameter table}
\begin{center}
\resizebox{\columnwidth}{!}{\begin{tabular}{c c l|c c l}
\toprule
\textbf{Module} & \textbf{Parameter} & \textbf{Description} & \textbf{Module} & \textbf{Parameter} & \textbf{Description} \\
\hline
Module B & $S_{neu}$ & Sites of detected neurons in $\phi$ &
Module D & $\zeta_{tip}$ & Tip intensity threshold \\
{} & $S_{gc}$ & Site of growth cones in $\phi$ &
{} & $n_{tips}$ & Number of neurite tips detected \\ \cline{4-6}
{} & $n_{neu}$ & Number of neurons &
Module E & $l_{seg}$ & Evolving neurite segment length \\
{} & $S_{tips}$ & Sites of detected tips &
{} & $S_{tr}$ & Site of neurite tracings \\  
{} & $l_{gc}$ & Size of growth cone &
{} & $l_{neu}$ & Neurite length \\  
{} & $\mathcal{H}_{conv}$ & Convolution box filter &
{} & $d_{geo}$ & Calculated geodesic distance \\  \cline{1-3}
Module C & $P_{initial}$ & Initial neuron cell center &
{} & $\theta_{t}$ & Calculated relative turning angle \\
{} & $ l_{total}$ & Evolving total neurite length &
{} & $\tau$ & Calculated tortuosity \\ 
{} & $\bar l_{total}$ & Average total neurite length &
{} & ($x_{tip}, y_{tip}$) & Neurite tip coordinates \\  
{} & $\bar \theta_{t}$ & Average relative turning angle &
{} & ($x_{root}, y_{root}$) &  Neurite root coordinates \\ \cline{4-6}
{} & $\bar \tau$ & Average tortuosity &
Algorithms & $\zeta_{soma}$ & Soma geodesic threshold \\
{} & $\bar n_{e}$ & Average number of neurite tips &
~\ref{Algorithm: Neurite_tracing_algorithm}\&\ref{Algorithm: tip_detection} & $I$ & Intensity field \\
{} & $\bar l_{seg}$ & Average neurite segment length &
{} & $A_{tip}$ & Area of the tip \\  
{} & {} & {} &
{} & $\gamma_{tip}$ & Tip area threshold \\
\bottomrule
\\[-1em] 
\multicolumn{6}{l}{Note: Parameters ($\bar l_{total}$, $\bar \theta_t$, $\bar \tau$, $\bar n_e$, and $\bar l_{seg}$) with a bar on top are obtained from experimental neuron cultures.} \\
\end{tabular}}
\end{center}
\vspace{-0.8cm}
\end{table}

\begin{figure}[ht]
\centering
\includegraphics[width=\linewidth]{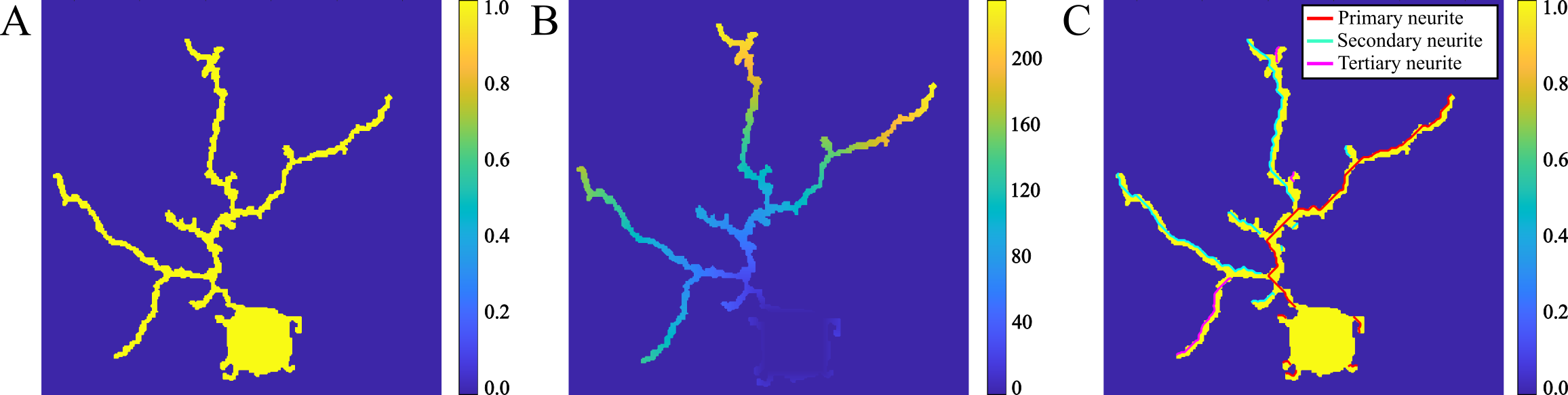}
\caption{Neurite tracing procedure (Algorithm~\ref{Algorithm: Neurite_tracing_algorithm}) that takes $\phi$ as input and traces three generations of neurites. (A) Input phase field $\phi$. (B) Geodesic distance $d_{geo}$ in neuron from $P_{initial}$ to neurite tips $S_{tips}$. Tracing is achieved by propagating along neurites based on $d_{geo}$. (C) Output all neurite tracings $S_{tr}$ and the sum of all neurite lengths.}
\label{figure: Neurite_tracing_procedure}
\end{figure}

\textbf{Growth cone calculation (Module B)}. Based on the IGA-C neuron growth solver (Module A), we incorporate the effect of various neurite morphometric features into the model by adjusting competitive tubulin consumption $\frac{dL}{dt}$ at growth cones of neurite tips. This is achieved by adjusting $M_{\phi}$, $r_g$, and $s_g$ values to change tubulin concentration balance at neurite growth cones $S_{gc}$. $S_{gc}$ is determined using a feature-driven neuron growth procedure (Modules B-D) that takes $\phi$ as input and outputs corresponding $S_{gc}$ for each neuron based on neurite morphometric features of the corresponding DIV in Table~\ref{table: neurite_feature_data}. During the neuron growth simulation, we run the IGA-C neuron growth solver until the final DIV is reached by iteratively determining the evolving DIV based on the total neurite length $\bar l_{total}$ (Figure~\ref{figure: Neurite_tracing_procedure}). In each iteration, we apply the feature-driven neuron growth procedure to adjust $\frac{dL}{dt}$, which affects the magnitude of $E$ (Eq.~\ref{energy equation}) for the phase field governing equation (Eq.~\ref{phase field equation}) to drive neurite outgrowth. During the procedure, we first use MATLAB's connected component analysis to identify neurons $S_{neu}$ within the phase field $\phi$. Then, for each identified neuron, we search for potential neurite tips and identify their locations $S_{tips}$ using a convolution operation with a box filter $\mathcal{H}_{conv}$, which will be explained later in Module D. With identified neurite tips $S_{tips}$, we can generate growth cones $S_{gc}$ by selecting $S_{tips}$ neighboring area based on the growth cone size $l_{gc}$ (Module E). To simulate axon differentiation, we adjust mobility term $M_{\phi}$ on $S_{gc}$ furthest from the initial neural cell center $P_{initial}$ to $100$ while keeping $M_{\phi}$ on the rest $S_{gc}$ as default $60$. The $M_{\phi}$ magnitude difference among neurite tips allows the longest neurite to grow faster than the rest, therefore achieving axon differentiation. We then increase the $r_g$ value to $50$ and decrease the $s_g$ value to $0$ to drive neurite outgrowth at $S_{gc}$ as shown in Module B. Based on semi-automated neurite morphometric features \cite{liao_quantitative_2022,liao_semi-automated_2022} shown in Table~\ref{table: neurite_feature_data}, we can obtain relative turning angle $\bar \theta_t$, average tortuosity $\bar \tau$, neurite endpoints $\bar n_e$, average segment length $\bar l_{seg}$, and the total neurite length, $\bar l_{total}$ to drive and constrain the neuron growth model. 

\textbf{DIV determination (Module C)}. To achieve intrinsic DIV transition during simulation, we need to correctly determine the evolving DIV for $S_{neu}$ during the simulation. As shown in Module C, we incorporate an automatic neurite tracing algorithm (Algorithm \ref{Algorithm: Neurite_tracing_algorithm}) into the computational pipeline to calculate the total neurite length, $l_{total}$, and use it to determine the current DIV. Algorithm \ref{Algorithm: Neurite_tracing_algorithm} takes neuron growth phase field $\phi$ as input and outputs all primary, secondary, and tertiary neurite tracings based on their morphology connections in the neurite structure. In this process, we first calculate the corresponding neuron cell initial coordinate $P_{initial}$ based on centroids of detected $S_{neu}$. Then, the algorithm calculates the geodesic distance, $d_{geo}$, within the neuron cell from $P_{initial}$. To correctly differentiate neurites from the soma in a binary domain, we leverage a soma geodesic threshold $\zeta_{soma}$. Neurite tracing is then achieved by propagating from the maximum geodesic value, $max(d_{geo})$, to the minimum geodesic value, $max(d_{geo})$, in the same connected neurite. As shown in Figure~\ref{figure: Neurite_tracing_procedure}, once primary neurite tracing is complete, the algorithm repeats the tracing procedure and removes traced neurites to differentiate among primary, secondary, and tertiary tracings. Finally, based on obtained neurite tracings, we can calculate the total neurite length, $l_{total}$, by summing up the geodesic distance of each neurite, $l_{neu}$, and determine the corresponding DIV.

 \algrenewcommand\algorithmicrequire{\textbf{Procedure}}
\begin{algorithm}[ht]
\textbf{Input}: Neurite growth field $\phi$\\
\textbf{Output}: Neurite tracing and total neurite length $l_{total}$
\begin{algorithmic}[1]
\caption{Automatic Neurite Tracing (Figure~\ref{figure: Neurite_tracing_procedure})} \label{Algorithm: Neurite_tracing_algorithm}
\Require{Primary neurite tracing}
\State Calculate the geodesic distance, $d_{geo}$, from the center of soma $P_{initial}$ to each tip
\State Remove soma region: Set $\phi = 0$ where $d_{geo}>\zeta_{soma}$
\For{each neurite}
\State Propagate and trace from $max(d_{geo})$ to $min(d_{geo})$ to obtain neurite tracings $S_{tr}$
\EndFor
\State Set $d_{geo}$ along $S_{tr}$ to 0 to remove traced neurites
\Require{Secondary and tertiary neurite tracing}
\State Repeat Step $3\&4$ to obtain secondary and tertiary neurites
\Require{Total neurite length calculation}
\For{each neurite tracing}
\State Calculate each neurite length: $l_{neu} = max(d_{geo}) - min(d_{geo})$
\EndFor
\State Compute total neurite length: $l_{total} = \sum{l_{neu}}$
\end{algorithmic}
\end{algorithm}

\algrenewcommand\algorithmicrequire{\textbf{Procedure}}
\begin{algorithm}[!hb]
\textbf{Input}: Neurite growth phase field $\phi$, $\zeta_{tip}=0.4$, and $\gamma_{tip}= 10$\\
\textbf{Output}: Detected tips $S_{tips}$
\begin{algorithmic}[1]
\caption{Neurite tip detection (Figure~\ref{figure: Tip_detection})} \label{Algorithm: tip_detection}
\Require{Detect potential neurite tip locations $S_{tips}$}
\State \textbf{Goal}: To obtain tips by extracting areas with intensity $I$ lower than $\zeta_{tip}$
\State Calculate intensity $I$ using convolution on $\phi$: $I = \phi \otimes \mathcal{H}_{conv}$
\If{$I > \zeta_{tip}\times\max(I)$}
    \State $I = 0$
\EndIf
\Require{Filter out tips with small areas (noise)}
\State \textbf{Goal}: To obtain tip centroids $S_{tips}$ and tip areas $A_{tip}$ using connected component analysis
    \For{each detected tip}
        \If{$A_{tip} \le \gamma_{tip}$}
            \State $S_{tips} = 0$
        \EndIf
    \EndFor
\end{algorithmic}
\end{algorithm}

\begin{figure}[!ht]
\centering
\includegraphics[width=\linewidth]{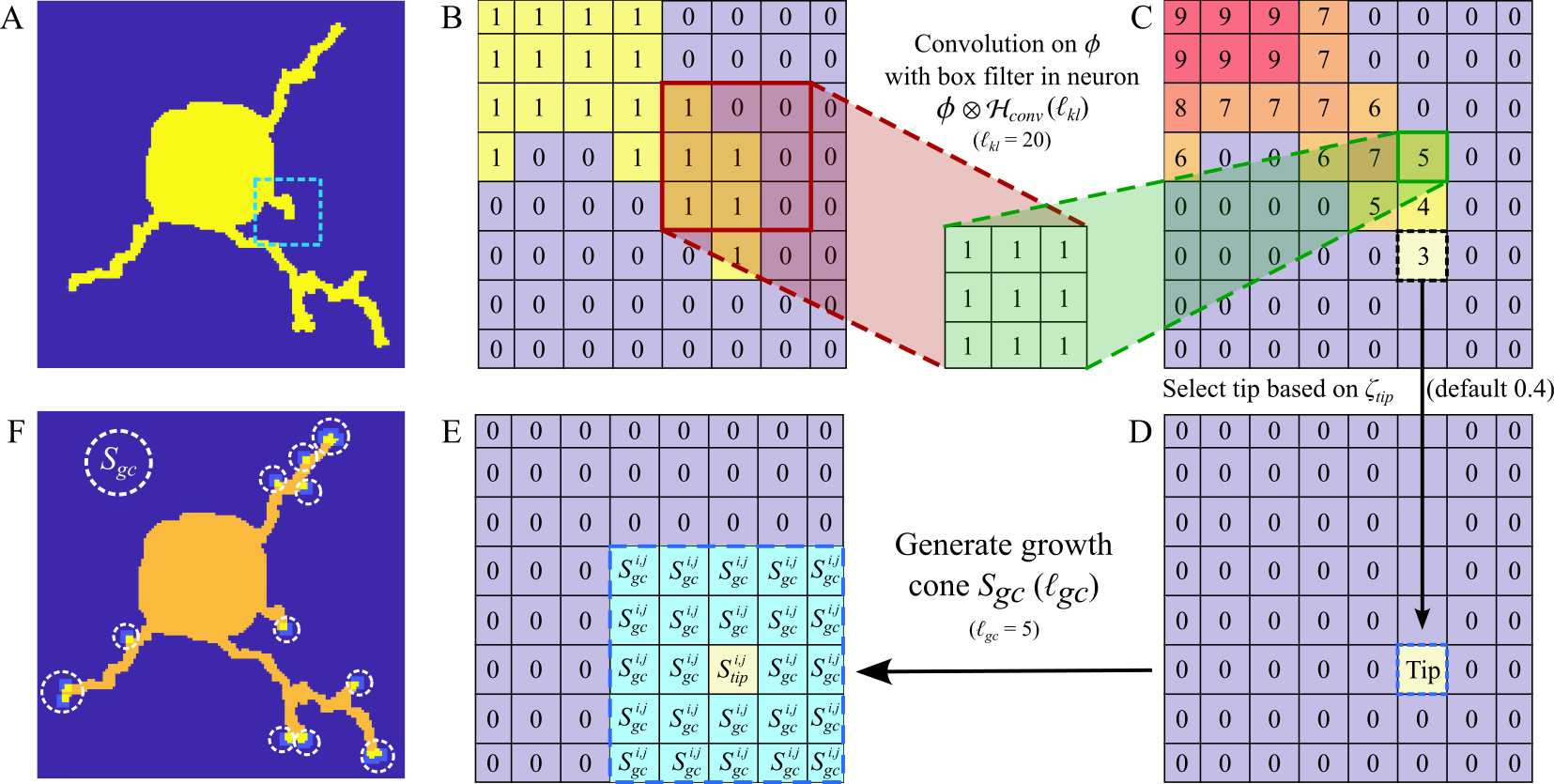}
\caption{Tip detection procedure (Algorithm \ref{Algorithm: tip_detection}). (A) Input $\phi$ field. (B) A zoomed-in view of the tip area (the cyan dashed square in (A)). (C) Intensity field $I$ is calculated using convolution to $\phi$ inside the neuron cell (where $\phi=1$). For visualization, $l_{kl}$ is set as 3 (default value is 20). (D) Detecting potential neurite tips $S_{tips}$ based on tip intensity threshold $\zeta_{tip} \times max(I)$. $\zeta_{tip}$ is set to $0.4$. (E) Growth cone $S_{gc}$ generated based on $S_{tips}$. $l_{gc}$ is set as 5. (F) Generated growth cone $S_{gc}$ across the entire neuron.}
\label{figure: Tip_detection}
\end{figure}

\textbf{Potential tip detection (Module D)}. To simulate neurite outgrowth in the IGA-C phase field neuron growth model, we carry out a ``growth-cone" like activation of the energy term $E$, as shown in Module D. This approach allows us to incorporate the number of neurite endpoints $\bar n_e$ from experiments into the neuron growth model. We develop a neurite tip detection algorithm (Algorithm \ref{Algorithm: tip_detection}) with an adjustable threshold value $\zeta_{tip}$ to constrain the number of tips $S_{tips}$. The tip detection algorithm can automatically detect neurite tips based on a given $\zeta_{tip}$ of $0.4$. In Figure~\ref{figure: Tip_detection}, the algorithm takes neurite growth pattern $\phi$ as input and outputs a list of potential tips in the domain. The neurite tips are detected as the locations of high intensity magnitude with fewer neighbors than the rest of the neurite regions (Figure~\ref{figure: Tip_detection}C). First, the algorithm calculates an intensity field $I$ by applying convolution on $\phi$ using a box filter $\mathcal{H}_{conv}$:
\begin{align}
    I(m,n) = \phi (m,n) \otimes \mathcal{H}_{conv}(m,n) =  \sum_{i=1}^{l_{kl}} \sum_{j=1}^{l_{kl}} \mathcal{H}_{conv}(i,j) \cdot \phi (m-i,n-j), \qquad 1\le m \le M, \quad 1\le n \le N,
    \label{convolution equation}
\end{align} where $\mathcal{H}_{conv}$ is the convolution kernel, and $l_{kl}$ is kernel size. $M$ and $N$ are the dimensions fo the domain. Then, we can identify tips by extracting sites with intensity magnitudes lower than $\zeta_{tip}\times max(I)$ (Figure~\ref{figure: Tip_detection}D). Based on the $1^{st}$ and $3^{rd}$ quartiles of $\bar n_e$ given in Table~\ref{table: neurite_feature_data}, we can adjust $\zeta_{tip}$ to constrain the number of detected tips to fit within the range. Because the phase field method solves interface movement based on energy balance, there exist small protrusions that could be incorrectly identified as tips by Algorithm~\ref{Algorithm: tip_detection}. Therefore, we apply a clearing step at the end to filter out tips with small areas based on a tip area threshold $\gamma_{tip}$ (the default value is 10). Furthermore, because Algorithm~\ref{Algorithm: tip_detection} is intensity-based with connected component analysis, neurites from the same neuron will not intersect as $I$ will decrease when neurites approach each other, while $I$ on neurites from different neurons will not be affected.

\begin{figure}[!ht]
\centering
\includegraphics[width=\linewidth]{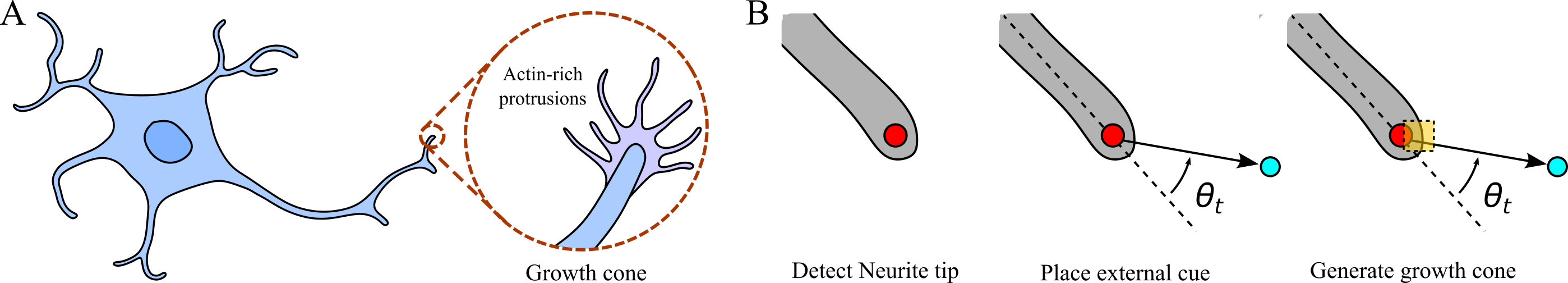}
\caption{Neurite growth cone and external cue-guided mechanism. (A) The growth cone is an actin-rich protruding area at the neurite tip. (B) External cue-guided mechanism. The neurite tip is first detected, then an external cue (cyan dot) is placed based on $\theta_t$, and finally a growth cone (orange square zone) is selected close to the external cue.}
\label{figure: External_cue_guided_mechanism}
\end{figure}

\textbf{Feature-based tip adjustment (Module E)}. After determining evolving DIV and detecting potential neurite tips, we use the external cue-guided mechanism and energy activation zones introduced in the neuron growth model \cite{qian_modeling_2022} to guide neurite growth. During neuron growth, an actin-rich protrusion at the neurite tip called the growth cone is responsible for exploring the surrounding environment and guiding the neurite growth toward extracellular cues (Figure~\ref{figure: External_cue_guided_mechanism}A). The extracellular medium will affect neurite growth because neurites consume various proteins and chemical concentrations during the growth and tend to grow towards extracellular attractive cues with high concentration. With the external cue-guided mechanism, we can specify the direction neurites should grow towards by changing the external cue placement (Figure~\ref{figure: External_cue_guided_mechanism}B). Based on the evolving DIV, we first select and implement the corresponding neurite morphometric features ($\bar \theta_t, \bar \tau,\bar l_{seg}$) from Table~\ref{table: neurite_feature_data}. During the feature-based tip guidance procedure, we can selectively update external cues to guide the growing neurites so that neurite segments between updates satisfy the range set by the $1^{st}$ and $3^{rd}$ quartiles ($Q_1$, $Q_3$) of $\bar l_{seg}$ by monitoring the average segment length $l_{seg}$. $l_{seg}$ was calculated by dividing $l_{total}$ by the number of neurite tips recorded. To update external cue placement in Module E, we identify neurite tracing $S_{tr}$ based on the neurite tip $S_{tips}$. Then, we calculate the neurite tortuosity $\tau$ by dividing the geodesic distances $l_{neu}$ between the tip and the root by the corresponding Euclidean distance:
\begin{align}
    l_{neu} = max(d_{geo})-min(d_{geo}),\\
    \tau = \frac{l_{neu}}{\sqrt{(x_{tip}-x_{root})^2+(y_{tip}-y_{root})^2}},
\end{align} where ($x_{tip},y_{tip}$) and ($x_{root},y_{root}$) are the coordinates of the neurite tip and root, respectively. Then, to implement the average experimental neurite tortuosity $\bar \tau$, we repeatedly generate new $\theta_t$ based on $\bar \theta_t$ to place the external cue, until the calculated neurite tortuosity $\tau \in [\bar \tau(Q_{1}),\bar \tau(Q_{3})]$. To generate a new $\theta_t$, we first calculate the direction of the neurite using endpoints of $S_{tr}$. Then, we generate a new $\theta_t$ using Gaussian distribution for each $S_{tips}$ based on mean $\mu$ and standard deviation $\sigma$ of experiment $\bar \theta_t$:
\begin{align}
  \theta_t(x) = \frac{1}{\sigma\sqrt{2\pi}} 
  \exp\left( -\frac{1}{2}\left(\frac{x-\mu}{\sigma}\right)^{\!2}\,\right).
\end{align} As shown in Figure~\ref{figure: External_cue_guided_mechanism}B, the external cue (cyan dot) is then placed near the neurite tip (red dot) based on $\theta_t$, and the tip section $S_{tip}$ closest to the cue is selected as the energy activation zone to simulate the effect of the growth cone. Once all neurite tips are detected, the updated $S_{tips}$ of each neuron are then passed back for growth cone generation in Module B.

\subsection{Neurite features driven neuron growth results}

In this section, we study single- and multiple-neuron growth scenarios using the proposed feature-driven neuron growth model and summarize improvements based on a comparison with experimental neuron cultures and our previous neuron growth model \cite{qian_modeling_2022}. Our feature-driven neuron growth IGA-C solver is developed using MATLAB (R2022a) on a desktop with 12 cores and 64GB memory. The batch computation used Bridges-2 Supercomputer at Pittsburgh Supercomputer Center through Advanced Cyberinfrastructure Coordination Ecosystem (ACCESS) \cite{xsede,ecss}. 

\begin{figure}[htbp]
\centering
\includegraphics[width=\linewidth]{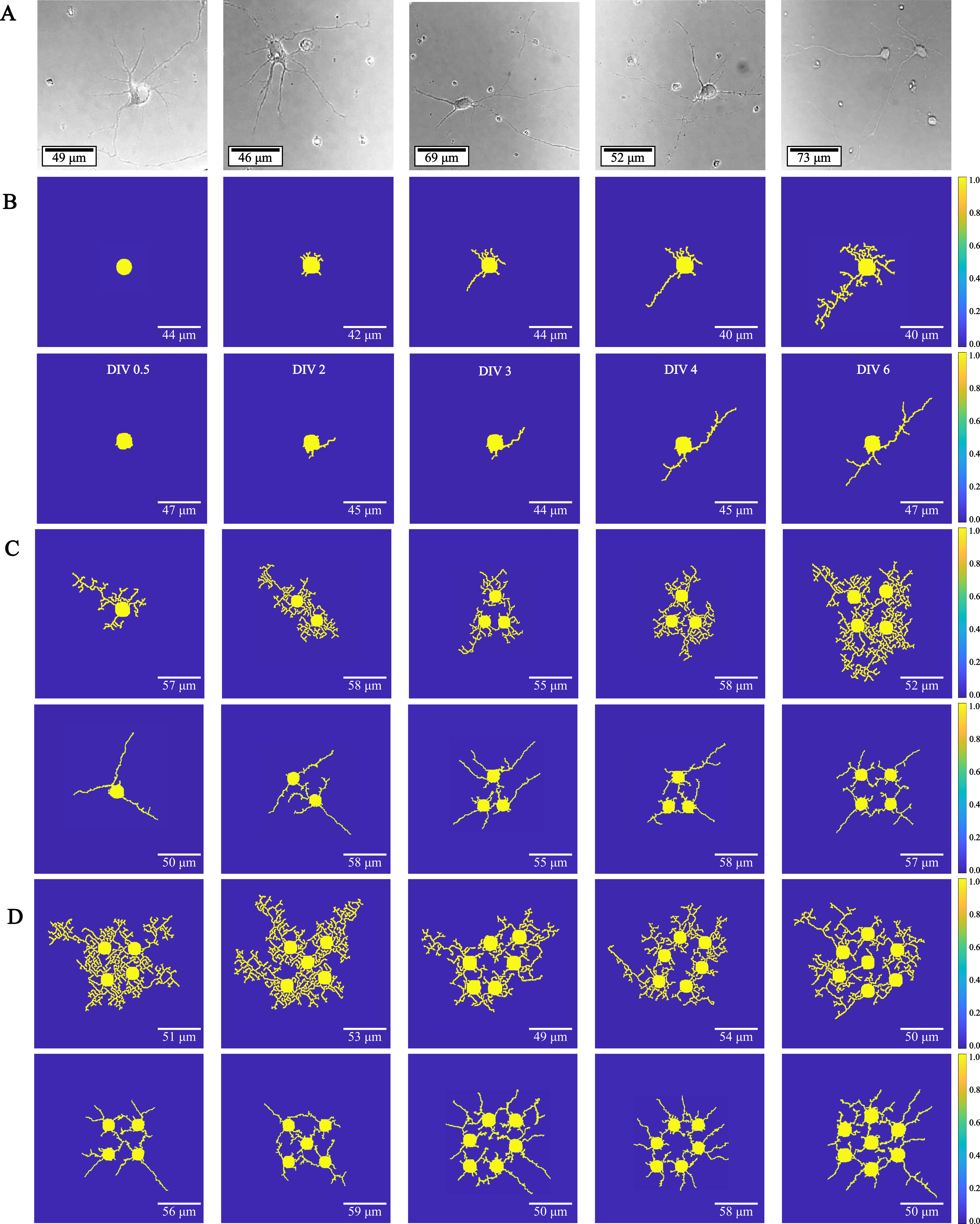}
\caption{Neurite morphometric features-driven neuron growth results with intrinsic growth stage transition. (A) Experimental neuron growth culture images of rat hippocampal neurons. (B) Neuron growth simulation results at each DIV using our previous neuron growth model \cite{qian_modeling_2022} (top row) and feature-driven neuron growth model (bottom row). (C-D) Comparison of both models for 1 to 7-neuron scenarios.}
\label{figure: neurite_feature_results}
\end{figure}

\begin{figure}[!ht]
\centering
\includegraphics[width=\linewidth]{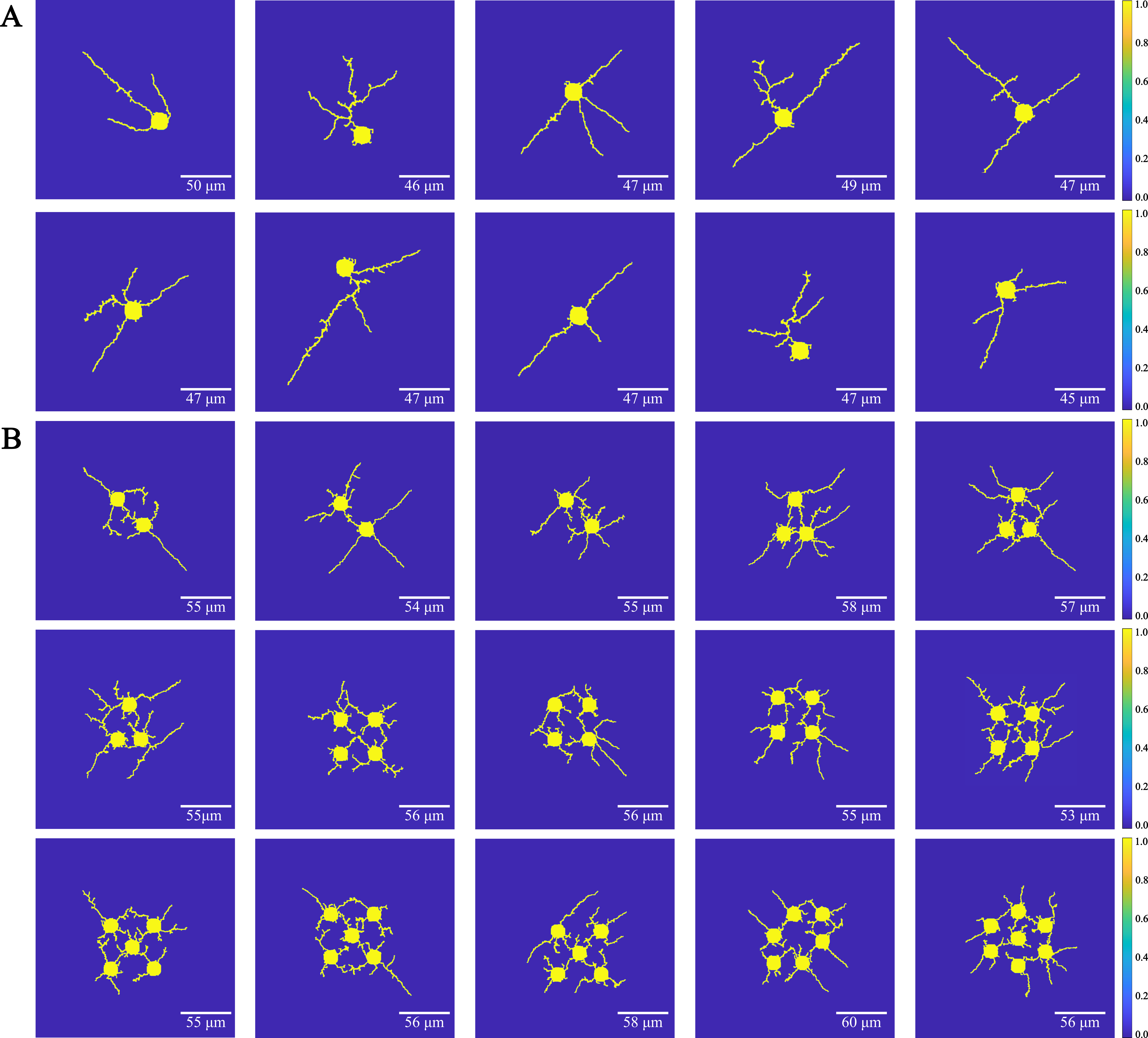}
\caption{Feature-driven neuron growth model simulation results. (A) DIV 6 single neuron growth simulation results with neurite morphometric features implementations. (B) Multiple neuron scenarios (2-7 neurons) simulation results.}
\label{figure: neurite_features_results_0}
\vspace{-0.3cm}
\end{figure} 

To compare the feature-driven IGA-C model with our previous neuron growth model, we set up scenarios for 1-7 neurons.  In Figure~\ref{figure: neurite_feature_results}A, we first show a selection of neuron culture images with different numbers of neurons as a reference. In Figure~\ref{figure: neurite_feature_results}B-D, we compare our previous neuron growth model (top row) and our feature-driven model  (bottom row) on a single neuron over the growth process from DIV 0.5 to DIV 6. It is obvious that our previous neuron growth shows unnecessary branching at a later stage not observed in experimental images (Figure~\ref{figure: neurite_feature_results}A), while our feature-driven neuron growth model results exhibit less neurite branching with much fewer neurite endpoints. This behavior illustrates the effect of  $\bar n_e$, which constrains the number of neurite endpoints throughout the growth. Our feature-driven model shows relatively straight neurites through each DIV. This indicates the effect of tortuosity $\bar \tau$, which continuously guides the neurite to grow in a relatively straight path through the external-cue-guided mechanism. To simulate multiple neurons (2-7 neurons) as shown in Figure~\ref{figure: neurite_feature_results}C\&D, we leverage MATLAB's connected component analysis to iterate through each neuron $S_{neu}$ within $\phi$ and apply Modules B-E accordingly. For each $S_{neu}$, individual $S_{tips}$ and $M_{\phi}$ are calculated and then combined together with other neurons at the end of the iteration to drive $\phi$ interface change in Module A.

In Figure~\ref{figure: neurite_features_results_0}, we present more results using the feature-driven neuron growth model. All neuron growth simulation results ranging from 1 to 7 neurons show similar morphology compared to corresponding experimental neuron cultures. Among single neuron scenarios (Figure~\ref{figure: neurite_features_results_0}A), most neurites attain straight neurite behaviors and experience oscillating small bending behaviors due to $\bar \tau$ implementations. Neurite branching behaviors are also observed as neurites grow further out. Some cases develop all three generations of neurites (primary, secondary, and tertiary). Both of these features are representative of increasing $\bar n_e$ at later DIV.

\section{CNN-based prediction of the neuron growth process}

In this section, we demonstrate a convolutional neural network (CNN) based surrogate model to predict neuron growth scenarios. We first introduce our CNN auto-encoder implementation with multiple layers. Then, we explain our data generation procedures. Finally, we present our model prediction results compared to the IGA-C solver.

\subsection{CNN auto-encoder model}
\begin{figure}[ht]
\centering
\includegraphics[width=\linewidth]{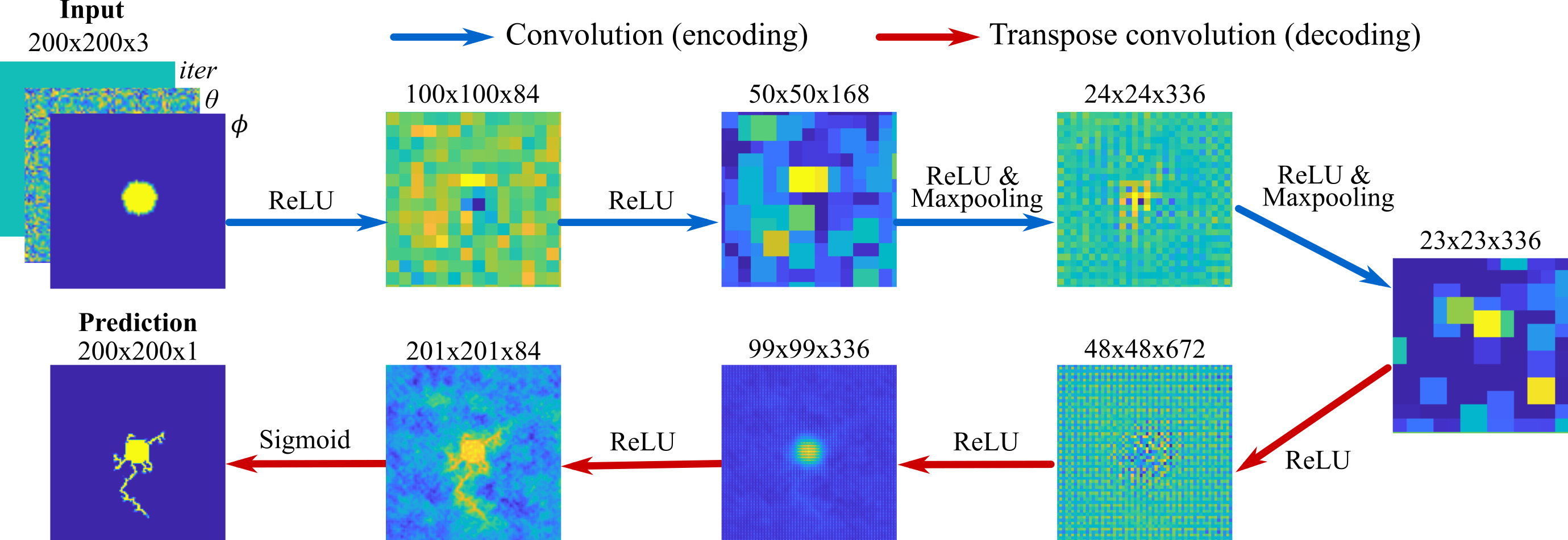}
\caption{Convolutional neural network (CNN) with the auto-encoder architecture for neuron growth prediction.}
\label{figure: CNN_model}
\end{figure}

We propose to use CNN with auto-encoder architecture to predict neuron growth. CNN is a subset of machine learning techniques proposed to operate on tensor data based on convolution operations. To achieve efficient and accurate prediction of the neurite growth pattern, we adopt a light convolutional autoencoder \cite{li_reaction_2020} as the backbone of our surrogate model, consisting of a multi-layer encoding and decoding architecture. The 3-channel input tensor consisting of phase field $\phi$, orientation $\theta$, and target iteration $iter$ goes through a 3-layer encoding process and then a 4-layer decoding process with a Sigmoid activation function at the output to predict neuron growth (Figure~\ref{figure: CNN_model}). 

Given a specific iteration, the model takes the neuron initializations as input and outputs the pattern of the grown neurite based on the target iteration number in the input tensors. The input includes the initial phase field $\phi$, the orientation $\theta$, and the target iteration number $iter$, where $\phi$ is a $300\times300$ matrix representing the initialized pattern of the neuron cell, and $\theta$ is randomly initialized between [0,1] for each neuron growth case. The model builds correlations among the input channels by extracting cross-channel features using convolution layers and rectified linear unit (ReLU) activation functions. The encoder reduces the spatial dimensionality of the input tensor and increases the number of channels (Figure~\ref{figure: CNN_model}). The decoder correlates and upsamples abstract features from the low-dimension representation to generate a $300\times300$ neuron growth prediction ($\phi$). During the encoding process, we also implement max-pooling layers to help prevent overfitting by only considering the element with the maximum value. As for the activation function, we choose ReLU because ReLU has shown similar accuracy while requiring less computational time compared to other widely used activation functions \cite{li_reaction_2020}. ReLU pushes negative values within the tensor to 0 while not affecting positive values. Because the phase field $\phi$ is binary, we add the Sigmoid activation function at the end to scale and reduce CNN model prediction to $0$ and $1$.

During the training process, we use binary cross-entropy (BCE) to calculate the loss gradients for back-propagation. Because the phase field $\phi$ is binary, we can treat the neuron growth problem as a binary classification problem in a 2D matrix for simplicity. BCE is ideally suited for this study since BCE loss computes the degree of divergence between binary predictions and the ground truth. Although the commonly used mean square error (MSE) loss function allows the model to make similar neurite growth pattern predictions, it cannot correctly generate predictions with magnitudes between 0 and 1. BCE is defined as:
\vspace{-0.2cm}
\begin{equation}
    BCE=-\frac{1}{N}\sum_1^N{[ y_{i} \times log( P(y_i))+(1- y_i) \times \log \left(P(1- y_i)\right)]},
    \label{BCE equation}
\end{equation} where $y_i$ is the binary label, and $P(y_i)$ is the probability of the entity being $y_i$ for all $N$ entities.  We train the model with a learning rate $\alpha$ of 1e-4 for 450 epochs. BCE enables us to treat the problem as a binary classification problem but does not provide a suitable error representation of the 2D domain. Thus, we use the mean magnitude of the relative error (MRE)) on the test dataset to better evaluate model performance and visualize relative errors of the 2D binary field ($\phi$) \cite{jorgensen_when_2022}. The MRE is defined as:

\vspace{-0.1cm}
\begin{equation}
    MRE=\frac{\sqrt{\frac{1}{n}\sum_1^n{(\phi_{gt}-\phi_{pred})^2}}}{\phi_{gt}^{max}-\phi_{gt}^{min}}\times100\%,
    \label{MRE equation}
\end{equation} where $n$ is the number of pixels in one matrix, $\phi_{gt}$ is the ground truth, and $\phi_{pred}$ is the prediction. $\phi_{gt}^{max}$ and $\phi_{gt}^{min}$ are the maximum and minimum values of the ground truth. From the generated dataset, each neuron case has three input channels: $\phi$ with neuron cells initialized near the center of the domain, orientation $\theta$ initialized to obtain evolving growth patterns, and target iteration $iter$s.

\subsection{Data generation}
Using the IGA-C-based phase field neuron growth model with incorporated neurite features, we run a batch of neuron growth cases with different neuron placements. For each case, we initialize neuron cells as solid circles inside the domain and randomly initialize an orientation $\theta$. For multiple neuron growth cases with complex neurite network formations, we initialize multiple neuron cells in the domain for the solver to develop neurite interactions among neurons. We consider seven types of neuron scenarios in our model, from 1 to 7 neurons, for generating the training dataset. In each case, we initialize the domain as a 60$\times$60 quadrilateral control mesh with neuron cell initialization and a corresponding random variable $\theta$ of the same size to introduce different growth behaviors into the model. The domain adaptively expands as the neurite growth approaches the boundary. We run the solver until DIV 6 is reached for each case with a time step $\triangle t = 0.05$. Simulation results at each iteration are extracted and stored based on collocation points in the domain. Because we adaptively expand the domain during the simulation, the final mesh size varies case by case. We downsample the domain into $300\times300$ to fit the dataset into our convolutional neural network. Then, to lower memory requirements and remove redundant information, we extract 60 data frames per case. Our CNN model takes the initial phase field $\phi$, case-specific randomly initialized orientation $\theta$, and target iteration $iter$ as the input and outputs predicted $\phi$ at target $iter$. Thus, we extract and assemble these three variables from each data point to obtain a dataset of size $(105,60,300,300,3)$. Before training, the dataset is shuffled and split 75/25 as training and testing dataset, respectively.

\subsection{Single and multiple neuron CNN predictions}
In this section, we present predictions of single and multiple neuron cases using our CNN surrogate model. After setting up our model architecture using Pytorch \cite{NEURIPS2019_9015}, we use Tesla V100-16 on Bridge-2 Supercompter at Pittsburgh Supercomputer Center (PSC) to train the model on $75\%$ of randomized neuron dataset for 450 epochs. The model reaches an average MRE of $2.56\%$ using our previous neuron growth dataset and $2.23\%$ using the feature-driven neuron growth dataset.

\begin{figure}[htbp]
\centering
\includegraphics[width=\linewidth]{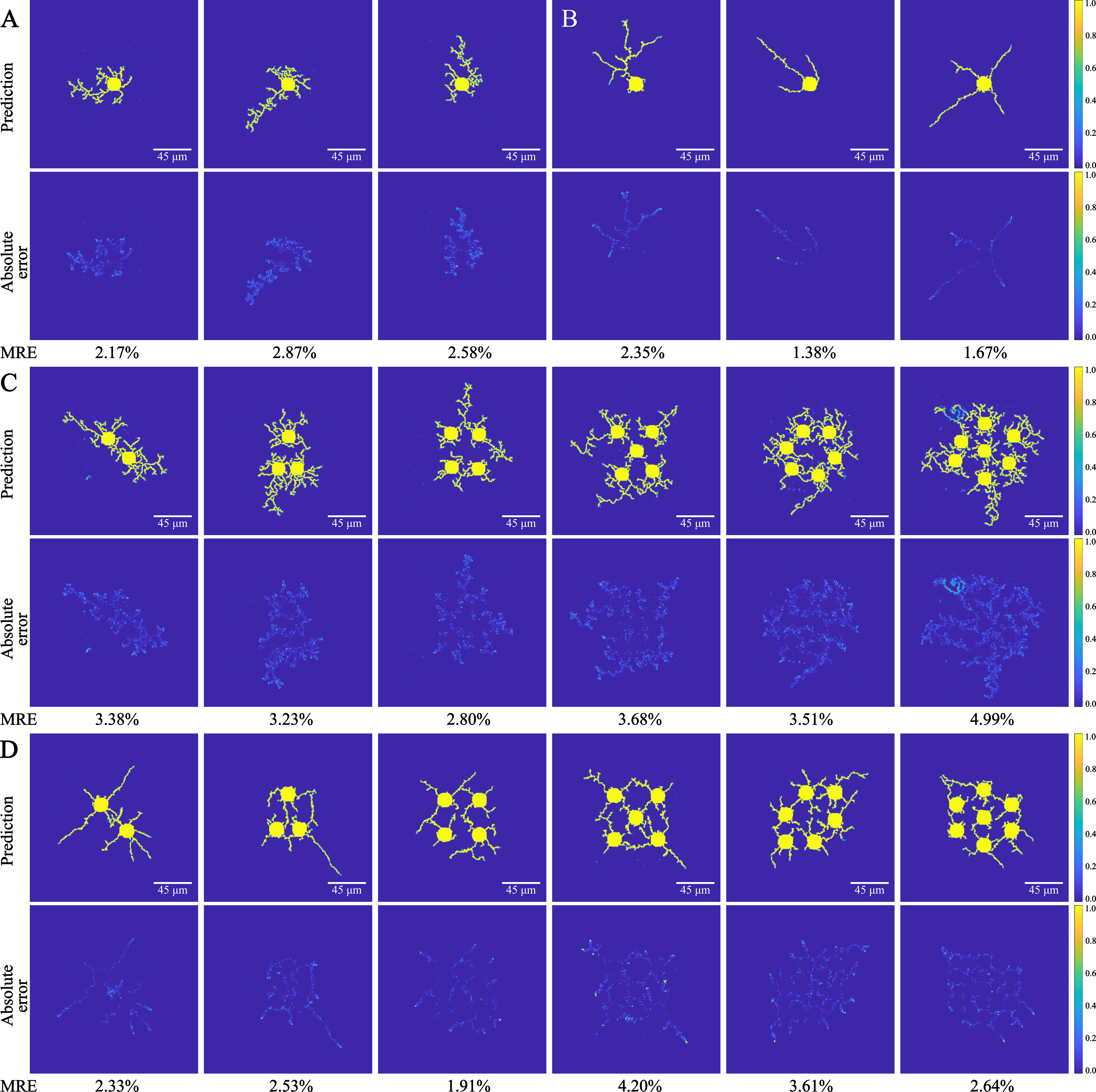}
\caption{CNN-based model prediction results of single neuron cases using two neuron growth models. (A, C) Single and multiple neuron growth prediction using our previous neuron growth model \cite{qian_modeling_2022}. (B, D) Single and multiple neuron growth prediction using feature-driven neuron growth model. The top row is the CNN prediction, and the bottom row is the corresponding absolute error. These results show that our CNN model can reproduce similar neurite growth patterns to the IGA-C solver with high accuracy (MRE $< 5\%$). }
\label{figure: CNN_results}
\end{figure}

Using our CNN-based model, we can predict the growth of single neuron cases and calculate the absolute error compared to corresponding IGA-C solver results (Figure~\ref{figure: CNN_results}A\&B). Each neuron prediction closely resembles neurite growth patterns shown in ground truth with complex branching and axon differentiation. We also observe that the highest error occurs at the neurite tips. As the $iter$ number increases, more neurite tips form due to branching. Consequently, the mean MRE of single neuron predictions goes up to $2.87\%$ in the cases shown. We then study multi-neuron scenarios (2-7 neurons) with complex neurite networks using our CNN-based model and compare the results with IGA-C results (Figure~\ref{figure: CNN_results}C\&D). In each case shown, we observe that more neurons form more complex neurite networks with the largest mean MRE of $4.99\%$ and that these results show that our model can reproduce similar neurite growth patterns to the IGA-C solver with high accuracy.

\begin{figure}[ht]
\centering
\includegraphics[width=\linewidth]{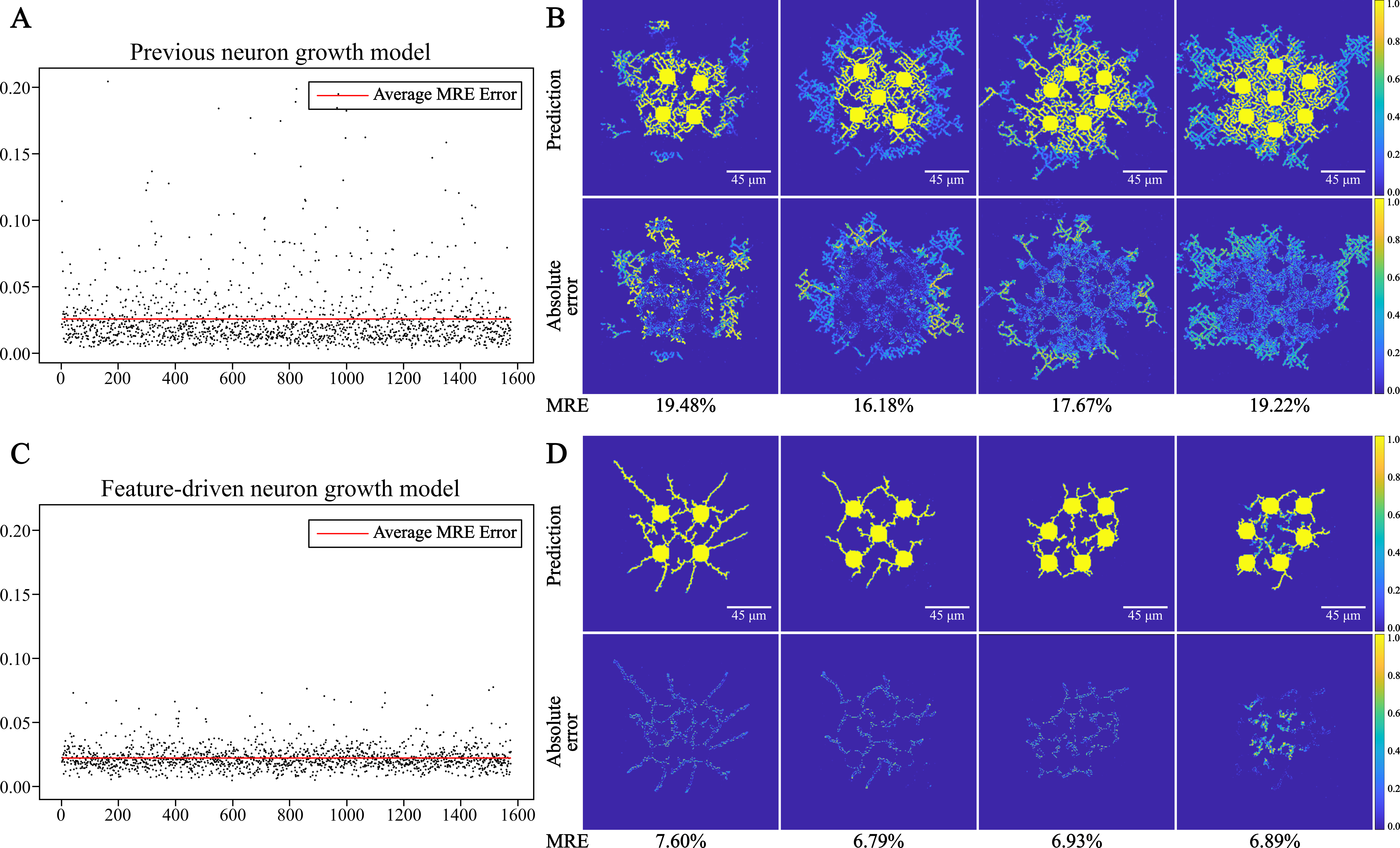}
\caption{CNN-based model prediction error. (A, C) CNN-based model accuracy statistics of our previous and feature-driven neuron growth models using the test dataset. MRE of all test datasets are plotted as scatter points, and the average MRE is plotted as the red line. (B, D) Four worst CNN predictions of both models.}
\label{figure: CNN_error}
\end{figure}

To better visualize our model performance, we plot the MRE in a scatter plot with the average MRE as a red line (see Figure~\ref{figure: CNN_error}A\&C). We observe that model prediction error climbs as $iter$ increases when complex neurite patterns form for all cases. The maximum MRE appears on multi-neuron cases at later $iter$, particularly for ones with complex neurite network formations. We record a maximum MRE of $19.48\%$ and $7.60\%$ for our previous neuron growth model and the proposed feature-driven neuron growth model, respectively. We select four worst-case predictions based on MRE for both models (see Figure~\ref{figure: CNN_error}B\&D) and observe that they exhibit similar incorrect neurite patterns in areas far away from the center, indicating a lack of information during training in these areas. This finding is within expectation because our neuron growth problem is set to grow from the center to the edge of the domain that is adaptively expanded during the simulation. Therefore, no matter how large our dataset is, there will always be a lack of information near the boundary. Consequently, our end-to-end CNN-based model may be limited when facing complex neurite networks near the domain boundary. We also observe that the computation time of our CNN prediction (approximately $0.25s$ per prediction) is vastly faster than that of the IGA-C neuron growth solver (approximately $240$ hours on average). Therefore, the CNN-based model can predict complex neuron growth patterns up to 7 orders of magnitude times faster than the IGA-C neuron growth solver. 

\section{Conclusion and future work}
This paper demonstrates a new computational pipeline to incorporate experimental neurite morphometric features into an IGA-C phase field model to simulate biomimetic intrinsic multi-stage neuron growth behaviors. Furthermore, we demonstrate that a CNN-based surrogate model can significantly reduce the associated computational cost for neuron growth predictions. Based on the results presented in this paper, we conclude that:
\begin{itemize}
    \item We have incorporated experimentally observed neurite morphometric features across 6 DIV into the IGA-C neuron growth model to drive and constrain the neuron growth process. The model generates biomimetic multiple-neuron scenarios with neurite interactions statistically comparable to experimental results. 
    \item The feature-driven neuron growth model is adaptable to different neuron growth behaviors by switching to a different set of neurite morphometric features. In this paper, we utilize rat hippocampal neuron data.
    \item  The CNN surrogate model can predict the multi-stage neuron growth process. Our model can learn from an abstract representation of neuron growth data and generate end-to-end accurate growth predictions at a given iteration. Our proposed model accurately predicts all growth stages ($<2.23\%$ error) while taking 7 orders of magnitude less computational times compared with our IGA-C neuron growth solver.
\end{itemize}

The CNN model significantly broadens future research possibilities due to its fast and accurate prediction capability. Our model enables researchers to visualize neuron growth results ahead of time without going through expensive neuron growth experimental trial and error to achieve desirable results. In particular, the proposed machine learning model will supplement traditional neuron culture experiments by enabling researchers to test and predict the effects of hypothetical experimental conditions before implementing them \textit{in vitro}. The proposed computational pipeline will allow researchers to plan experiments more efficiently and reduce experimental costs significantly by generating large amounts of predictions based on given initial conditions.

Whereas the model presented in this work substantially improves neuron growth modeling accuracy and speed, some limitations point to interesting future research directions. For example, the current features-driven IGA-C neuron growth model does not consider experimentally-measured growth rates of axons and dendrites \cite{dotti_establishment_1988}. Our growth model leverages the fixed mobility term difference to differentiate the growth rates of dendrites and axons. In the future, this can be improved by incorporating experimentally-measured growth rates. Another limitation is the lack of informative neuron growth factors, such as extracellular and intracellular concentrations, crucial for understanding the biophysics process behind the scenes. This information in future models may help reveal critical and previously unknown relationships between neuron growth and neurodegenerative diseases. In the future, we can combine our model with material transport in neurons to explore complex biophysics-coupled neurite morphologies \cite{li2022modeling_1,li2022modeling}. Finally, our feature-driven neuron growth model is limited to 2D using hierarchical B-splines. We plan to extend our neuron growth model to 3D and implement truncated hierarchical B-splines (THB-spline) \cite{pawar2018dthb3d_reg,wei2015truncated}. Note that our current neuron growth code is implemented in MATLAB, we plan to convert it into C/C++ to improve computational speed. Furthermore, the surrogate CNN model is trained based solely on simulation data without any information from the phase field governing equations. The residuals of the governing equations could be considered during the model training process using physics-informed loss functions \cite{raissi_physics-informed_2019,li_deep_2021}. A generalized surrogate model can be trained and adjusted to different types of neuron growth using the transfer learning approach \cite{goswami_transfer_2020}. Along the deep learning direction, we are also exploring more efficient and effective model architectures such as recurrent neural networks (RNN) and transformers \cite{NIPS2017_3f5ee243} to better understand temporal changes during the neuron growth process \cite{shi_convolutional_2015}.

\section*{Code and data availability}
The code and datasets generated and analyzed in this paper are accessible in the "FDNGCNN" GitHub repository. \url{https://github.com/CMU-CBML/FDNGCNN} (DOI:10.5281/zenodo.7853871). Correspondence and requests for code and data should be addressed to K.Q. or Y.J.Z.

\section*{Declaration of competing interest}
The authors declare no known competing financial interests or personal relationships that could have appeared to influence the work reported in this paper.

\section*{Acknowledgement}
K. Qian and Y. J. Zhang  were supported in part by the NSF Grant CMMI-1953323. All of the authors were supported in part by the PITA (Pennsylvania Infrastructure Technology Alliance) and PAMIP (Pennsylvania Manufacturing Innovation Program) grants. K. Qian was also supported by Bradford and Diane Smith Graduate Fellowship. In addition, A. S. Liao and V. A. Webster-Wood were supported in part by an NSF CAREER award ECCS-2044785. A. S. Liao was also supported by NSF Graduate Research Fellowship Grant DGE-1745016 and Carnegie Mellon University Jean-Francois and Catherine Heitz Scholarship. This work used RM-node and GPU-node on Bridges-2 Supercomputer at Pittsburgh Supercomputer Center \cite{xsede,ecss} through allocation ID eng170006p from the Advanced Cyberinfrastructure Coordination Ecosystem: Services \& Support (ACCESS) program, which is supported by National Science Foundation grants \#2138259, \#2138286, \#2138307, \#2137603, and \#2138296.

\bibliographystyle{elsarticle-num}
\bibliography{reference}
\end{document}